\documentstyle[12pt,leqno]{article}

\makeatletter
  \newcommand{\C}{{\bf C}}
  \newcommand{\Z}{{\bf Z}}
  \newcommand{\Met}{\mathop{\rm Met}_1}
  \newcommand{\R}{{\bf R}}
  \newcommand{\Hyp}{{\cal H}}
  
  \newcommand{\cmcone}{$\mathop{\rm CMC}$-$1$}
  \newcommand{\Fl}{\mathop{\it Fl}\nolimits}
  \newcommand{\fl}{\mathop{\it fl}\nolimits}
  \renewcommand{\sl}{\mathop{\it sl}}
  \newcommand{\PSL}{\mathop{\it PSL}}

  \newtheorem{Thm}{Theorem}
    
  \newtheorem{Prop}[Thm]{Proposition}
  
  \newtheorem{Cor}[Thm]{Corollary}
  \newtheorem{Example}[Thm]{Example}
    \newenvironment{Exa}{\begin{Example} \rm}{\end{Example}}
  \newtheorem{Definition}[Thm]{Definition}
    
  \newtheorem{Observation}{Observation}
    \newenvironment{Obs}{\begin{Observation} \rm}{\end{Observation}}
  \newtheorem{Remark}{Remark}
   \newenvironment{Rmk}{\begin{Remark} \rm}{\end{Remark}}
   
  \newtheorem{Acknowledgements}{Acknowledgements}
    \newenvironment{Ack}{\begin{Acknowledgements} \rm}{\end{Acknowledgements}}
  
    \newtheorem{Problem}{Problem}
  \newenvironment{Prob}{\begin{Problem} \rm}{\end{Problem}}

  \newenvironment{Proof*}[1]{\vskip.3cm\par\noindent%
    \mbox{\it #1.  }}{%
    \hfill \mbox{(q.~e.~d.)}\vskip.5cm}

  \newcommand{\eqref}[1]{\mbox{\rm (\ref{#1})}}
  \@addtoreset{equation}{section}
  \@addtoreset{footnote}{section}
\def\section{\@startsection {section}{1}{\z@}{3.5ex plus 1ex minus 
 .2ex}{2.3ex plus .2ex}{\large\bf}}
\makeatother
\title{ 
  A new flux for mean curvature 1 surfaces in
  hyperbolic 3-space, and applications
}
\date{ }
\author{\normalsize \bf
   \begin{tabular}[t]{c}
    Wayne Rossman\raisebox{0.8ex}{\tiny\rm 1},\\
    {\footnotesize Kyushu Univ.}
   \end{tabular}
   \begin{tabular}[t]{c}
    Masaaki Umehara\raisebox{0.8ex}{\tiny\rm 1}\\
    {\footnotesize Osaka Univ.}
   \end{tabular}
   and 
   \begin{tabular}[t]{c}
    Kotaro Yamada\raisebox{0.8ex}{\tiny\rm 1,2}\\
    {\footnotesize Kumamoto Univ.}
   \end{tabular}
}
\begin{document}
\maketitle

\begin{abstract}
  Using the Bryant representation,  
  we define a new flux on homology classes of \cmcone{} 
  surfaces in $\Hyp^3$, satisfying a balancing formula
  which is useful to show nonexistence of certain kinds of 
  complete \cmcone{} surfaces.
\end{abstract}

\addtocounter{footnote}{1}
\footnotetext{Supported by Volkswagen-Stiftung (RiP Program in 
              Mathematisches Forschungsinstitut Oberwolfach).}
\addtocounter{footnote}{1}
\footnotetext{Supported by Inamori Foundation.}
\setcounter{footnote}{0}

\noindent 
{ \large \bf Introduction}\,\,\,
Recently, flux and balancing formulas have been found very
useful techniques for proving results about surfaces of 
constant mean curvature $H$ (CMC-$H$) in 3-dimensional space forms.  
For CMC surfaces in Euclidean space $\R^3$ 
a flux has been defined on homology classes 
of closed curves $\gamma$ within the surfaces \cite{kks}, 
which we call the KKS-flux here.  
Also, in \cite{kkms}, 
a corresponding flux is defined for CMC surfaces in 
hyperbolic space $\Hyp^3$ of constant sectional curvature $-1$, 
which we call the KKMS-flux.  

In this work we define a new flux $\Fl$ using 
the Bryant representation \cite{Bry}.  
$\Fl$ satisfies a balancing formula (Theorem~\ref{thm:1.1}). 
We give a condition for when a regular embedded end has nonzero $\Fl$ flux 
(Corollary~\ref{cor:1.5}),
which is useful to show nonexistence of complete \cmcone{} surfaces under
certain assumptions about the behavior of their ends 
(Corollary~\ref{cor:1.3} and Theorem~\ref{thm:1.8}).  
 
The need for this new definition of 
flux for \cmcone{} surfaces in $\Hyp^3$ becomes 
clear if one looks carefully at the KKS-flux and KKMS-flux.  
The KKS-flux (resp.~KKMS-flux) is the sum of an integral along
a curve $\gamma$ in the surface and an integral over 
a ``cap" which is some 
2-dimensional disk in $\R^3$ (resp.~$\Hyp^3$) bounded 
by $\gamma$.  When $H=0$, the integral over the cap automatically 
vanishes.
When $H \neq 0$ (resp.~$|H|>1$) 
and the surface is of finite topology in $\R^3$ (resp.~$\Hyp^3$), 
the properly embedded annular ends of the 
surface converge asymptotically to Delaunay surfaces \cite{kks} 
(resp.~hyperbolic Delaunay surfaces \cite{kkms}), hence the curves
in neighborhoods of the ends can be kept at a uniformly bounded length, 
and the caps 
bounded by these curves can be kept to a uniformly bounded area. 
For this reason, 
when $|H| \ge 0$ (resp.~$|H|>1$) the 
flux is computable.  

However for a CMC-$H$ surface 
in $\Hyp^3$ with $0<|H|\leq 1$,
the cap has arbitrarily large area at an end. 
This makes it difficult to 
use the KKMS-flux effectively, in particular,
on complete \cmcone{} surfaces in $\Hyp^3$.  

For \cmcone{} surfaces in $\Hyp^3$ 
the flux $\Fl$ has the advantage of being 
computationally simpler than the KKMS-flux. 
Moreover, KKS-flux can be considered as a limit of our flux.
(See Remark after Theorem 1.)
Unfortunately, due to the difficulty of computing 
the KKMS-flux, any possible relationship between 
the KKMS-flux and $\Fl$ has not yet been determined.
However, both fluxes are $\sl(2,\C)$-valued, and this and other shared 
properties (see Example \ref{kumamotoexample}) 
suggest the possibility of a relationship. 

\begin{Ack}
We wish to thank Shin Kato for fruitful discussions
during our stay at Oberwolfach Institute.
We also wish to thank L-F. Cheung for his 
helpful comments and encouragement.
We are deeply grateful to the RiP program at
Oberwolfach Institute for its kind hospitality
and exemplary working conditions.
\end{Ack}
\section{ A flux formula for CMC-1 immersions}

Let $M$ be a Riemann surface.
We introduce the flux $\Fl$ for any conformal branched \cmcone{} 
immersion $x:M \to \Hyp^3$.  
(A point where the first fundamental form vanishes is called
a {\it branch point} of a \cmcone{} surface. Branch points on \cmcone{} 
surfaces are all isolated.)  
Let $\widetilde M$ be the universal cover of $M$.  
Take Weierstrass data $(g,\omega)$ of $x$ and let 
$F:\widetilde M\to \PSL(2,\C)$ be the unique lift of $x$ with respect to 
$(g,\omega)$,
which is a null holomorphic map (cf. \cite[Theorem 1.4]{UY3}).  
Let $G$ and $Q$ be the hyperbolic Gauss map and the Hopf differential of 
$x$.

The inverse matrix function,
$F^{-1}:\widetilde M\to \PSL(2,\C)$ is also a null holomorphic map,
and a new ``dual" branched \cmcone{} immersion
$x^\#:\widetilde M\to \Hyp^3$
is determined by
$x^\#:=F^{-1} \cdot (F^{-1})^*$.  
As shown in \cite{UY4}, this 
duality operation
interchanges the roles of the hyperbolic Gauss map and the secondary 
Gauss map. The Weierstrass data of the dual surface is given by
$(G,\omega^\#)$, where
\begin{equation}\label{2.1}
  \omega^\#:=-Q/dG. 
\end{equation}
In particular, the following identity holds:
\begin{equation}\label{def:dualalpha}
  \alpha^\#:= (F^{-1})^{-1}\cdot d(F^{-1})
            = - dF\cdot F^{-1} 
            = \left(
                 \begin{array}{rr} 
                       G& -G^2 \\
                       1& -G\hphantom{^2}  
                 \end{array}
              \right) \omega^\#.
\end{equation}
The key point is that $G$ and $Q$ are 
single-valued on $M$, and thus 
$\alpha^\#$ is single-valued on $M$.
So we can define a well-defined $\sl(2,\C)$-valued flux $\Fl(\gamma)$ 
for any loop $\gamma$ on $M$ by
\begin{equation}
  \Fl(\gamma):= \frac{-1}{2 \pi i} \int_{\gamma}dF\cdot F^{-1}.
\end{equation}
Since $dF\cdot F^{-1}$ is holomorphic on $M$,
the flux $\Fl(\gamma)$ depends only on the homology class
represented by the loop $\gamma$.

Now we suppose that there exists a closed
 Riemann surface $\overline{M}$
and a finite number of points $\{p_1,\ldots,p_n\}$ 
such that $M$ is biholomorphic to
$\overline{M}\setminus \{p_1,\ldots,p_n\}$ and the first fundamental form
is complete at each $p_j$. Each $p_j$ is called an {\it end} of the branched 
surface.  
(These conditions are automatically satisfied
when $x$ is an immersion which is complete 
and of finite total curvature. )
Then we have a well-defined $\sl(2,\C)$-valued flux $\Fl_j$ at each end $p_j$
by
\begin{equation}
  \Fl_j:=\Fl(\gamma_j) \;\;\;
  \left(= \frac{-1}{2 \pi i} \int_{\gamma_j}dF\cdot F^{-1}\right),
\end{equation}
where $\gamma_j$ is a loop surrounding the end $p_j$.
The following balancing formula is 
easily proved by triangulating $\overline{M}$ and summing the 
results of the Cauchy residue theorem.

\begin{Thm}\label{thm:1.1}
Let $x:M\to \Hyp^3$ be a branched \cmcone{} immersion.
Then the total sum of flux at the ends vanishes; that is, we have the 
balancing formula
\begin{equation}
  \sum_{j=1}^n \Fl_j=0 \; .
\end{equation}
\end{Thm}

\begin{Rmk}
For a minimal immersion $x:\overline{M}\setminus \{p_1,...,p_n\}\to \R^3$
with the Weierstrass data $(g,\omega)$, the KKS-flux is given
by the residue of the meromorphic form 
$((1-g^2), \sqrt{-1}(1+g^2),2g)\omega$ at each end.
As seen in [UY2], the Weierstrass representation can be 
considered as a limit of the Bryant representation.
Using this, one can easily see that KKS-flux is the limit of
our flux $\Fl$.
\end{Rmk}

Next we consider one regular end of a \cmcone{} immersion.
Using a local complex coordinate $(U,z)$ around an end 
$z=0$, the Hopf differential $Q $ is expressed as
\begin{equation}\label{1.7}
  Q(z)=q(z)\,dz^2=
       \left(\frac{q_{-2}}{z^2}+\frac{q_{-1}}{z}+\cdots\right)\,dz^2.
\end{equation}
An end $z=0$ is called a regular end of {\it Type I\/} 
if $q_{-2}\ne 0$, and called
a regular end of {\it Type II\/} if $q_{-2}=0$.  
(Regular ends are those at which the hyperbolic Gauss map is
meromorphic. An end is regular if and only if the Hopf differential
$Q$ is either holomorphic or has a pole of order at most 2 
at the end. [Bry], [UY1].)

\medskip
\noindent
{\bf Regular ends of Type I}\,\,
We compute the flux at an end $x$ of Type I.
We set $z=0$ at the end.  
Then the Weierstrass data $(G,\omega^\#)$ of the dual end 
$x^{\#}$ can be normalized as
\begin{equation}\label{3.1}
  G(z)=z^l \hat G(z), \qquad \omega^\#=z^k \hat w^\#(z) dz,
\qquad (\hat w^\#(0)\ne 0, \; \hat G(0) \ne 0),
\end{equation}
where $l>0$ and $k$ are integers, and $\hat G(z)$ and 
$\hat w^{\#}(z)$ are holomorphic functions.  
Since the metric of $x^{\#}$ is also complete (see Yu [Y] or \cite{UY4}),
by \cite{UY1} we have
$  k+1<0$. 
Since $q_{-2}\ne 0$, (1) yields that
$k+l=-1$. 
So the diagonal components $\pm G\omega^\#$ of $\alpha^{\#} = -dF \cdot 
F^{-1}$ have simple poles at $z=0$, and hence: 

\begin{Prop}\label{prop:1.2}
  A regular end of Type I has a non-vanishing flux.
\end{Prop}
\begin{Cor}\label{cor:1.3}
  There is no branched $1$-ended \cmcone{} surface with Type I end. 
\end{Cor}
  
\noindent
{\bf Regular ends of Type II}\,\,
Next we compute the flux for \cmcone{} ends of Type II.
Since this is not as easy as the Type I case,
 we do it only for some special cases, including the case of
embedded ends.  
For the case of minimal surfaces in $\R^3$, any
regular end of Type II has vanishing flux, 
but this is not true for the \cmcone{} case in $\Hyp^3$ (see 
Example \ref{kumamotoexample}).  

We fix a regular end of Type II at $z=0$.  
The Weierstrass data $(G,\omega^\#)$ of the dual end 
$x^{\#}$ can again be normalized with integers $l$ ($l>0$) and $k$ as
\begin{equation}\label{4.1}
  G(z)=z^l\, \hat G(z), \quad 
  \omega^\#=  z^k\, \hat w^\#(z)\, dz, 
  \qquad 
  (\hat w^\#(0)\ne 0, \; \hat G(0) \ne 0).  
\end{equation}
Since $q_{-2}=0$, Lemma~2 and Proposition~4 in \cite{UY4} yield 
that the secondary Gauss map $g$ is of the form
$
  g(z)=z^l\, \hat g(z),
$
where $\hat g(z)$ is a holomorphic function 
satisfying $\hat g(0)\ne 0$.  
Since $g$ is then single-valued at $z=0$,
$x^\#$ is also single-valued.  
(When $g$ and $G$ are both 
  single-valued on $U$, 
  then so are the lifts $F$ and $F^{-1}$ of $x$ and $x^{\#}$.
  See section 3 of \cite{Sm} or Theorem 1.6 of \cite{UY3}.)  
Hence by Theorem 2.4 in \cite{UY1}, the log-term coefficient $\theta$ of
the following differential equation vanishes: 
\begin{equation}\label{star}
  X''-\frac{(z^k\,\hat w^\#(z))^\prime}{z^k\,\hat w^\#(z)}X'+q(z)\,X=0,
\end{equation}
where $Q=q(z)dz^2$.
We set
\[
  \hat w^{\#}(z) = w_0 + w_1 z + w_2 z^2+\dots \; \; \; , 
\] 
where $w_0 \neq 0$.  
By the expressions \eqref{1.7} and \eqref{2.1}, we have
\begin{equation}\label{1.19}
  w_0=-\frac{q^{}_{l+k-1}}{l\, \hat G(0)}, \; \; 
  w_1 = -\frac{q^{}_{l+k}+(l+1) w_0 \,\hat G^\prime(0)}{l\, \hat G(0)}.  
\end{equation}
Let $\gamma$ be a loop surrounding the origin.
Since $q_{-2}=0$, $G\omega^{\#}=-Q\cdot (G/dG)$ and $G^2\omega^{\#}=-
GQ\cdot (G/dG)$ are holomorphic at $z=0$.
So, we have
\begin{equation}\label{eq:flux2}
  \frac{-1}{2\pi i}   
  \int_{\gamma}
       \left(
       \begin{array}{cc}
                G & -G^2\\
                1 & -G\hphantom{^2} 
       \end{array}
       \right) \omega^{\#}
  = -  \left(
       \begin{array}{cc}
               0 & 0 \\
               w_m & 0 
       \end{array}
       \right) \; \; , 
\end{equation}
where
$m=-k-1$ is the multiplicity of the end (see \cite{UY1}).
It can be easily seen that equation \eqref{star}
has no log-term if and only if (See Appendix of [RUY2] 
for the most direct reference, or
see [Bie] for a more detailed reference.)
\[
  \begin{array}{llll}
    w_1&=& - q_{-1} w_0  \qquad & (m=1),    \\
    4w_2 &=& - q_{0}w_0 - 3 q_{-1} w_1 - q_{-1}{}^2 w_0 \qquad 
                                & (m=2).
 \end{array}
\]
By \eqref{1.19} and \eqref{eq:flux2} and the above conditions, 
we have: 
\begin{Thm}\label{thm:1.4}
  A regular Type II end of multiplicity $m (=1 \mbox{ or\/ } 2)$
  has non-vanishing flux if and only if 
\small
$$ 
\begin{array}{lll}
   &    q_{-1}\ne 0  \qquad & (m=1),    \\
   &    q_{-1} = 0, \; q_0 \ne 0 \qquad \mbox{  or }   \\
   &\hspace{2em}
    q_{-1} \ne 0, \; \; 4\left(1-\displaystyle\frac{\hat G^\prime(0)}
        {\hat G(0)}\right)q_{0}+q_{-1}^2 \ne 0 \qquad & (m=2).    
 \end{array}
$$
\normalsize
\end{Thm}

Since an end is embedded if and only if $m=1$, by 
Proposition~\ref{prop:1.2} and Theorem~\ref{thm:1.4}, we have:  

\begin{Cor}\label{prop:1.4}\label{cor:1.5}
  A regular embedded end  has vanishing flux if and only if
  the Hopf differential $Q$ is holomorphic at the end. 
\end{Cor}

\begin{Exa}\label{example:1.7}
  We note that the flux
  does not depend only on 
  the asymptotic behavior of the end even in the Type I case.  
  For example, the catenoid cousin has the following data defined on 
  $M = \C \setminus \{ 0 \}$:  
\[
  g=z^\mu,\quad G=z,\quad Q=\frac{1-\mu^2}{4}\frac{dz^2}{z^2}
  \qquad(\mu\in\R^+\setminus\{1\}).
\]
The flux at $z=0$ is then given by
\[
  \frac{-1}{2 \pi i} \int \frac{1-\mu^2}{4} 
    \left(
    \begin{array}{ll} z^{-1} & -1 \\
              z^{-2} & -z^{-1}
    \end{array}
    \right)\,dz
    =\frac{\mu^2-1}{4}
    \left(
    \begin{array}{rr}
              1 & 0 \\
              0 & -1
    \end{array}
    \right)\;.
\]
Now we consider another example given by
\[
  g=z^\mu, \quad G=z+\frac{1}2z^2, \quad Q=\frac12(S(g)-S(G))
   =\frac{1-\mu^2}{4}\frac{1}{z^2}-\frac34\frac{1}{(1+z)^2}.
\]
($S$ is the Schwarzian derivative \cite{Bry}, \cite{UY1}.)  
This example is a branched surface (cf.\ \cite[Theorem 2.2]{UY3}).  
Since $G(0)=0$, 
one can check using (5.16) of \cite{UY1} that the end at the origin
is asymptotic to the end at $z=0$ of 
the above catenoid cousin.
The flux at the origin is
\[
  \frac{-1}{2 \pi i} \int
  \frac{1-\mu^2}{4}
  \left(
  \begin{array}{rr}
         \hphantom{-}z^{-1} & 0\hphantom{^{-1}} \\
         -z^{-1} & -z^{-1}  
  \end{array}
  \right)
  =
  \frac{\mu^2-1}{4}
  \left(
  \begin{array}{rr}
  1 & 0 \\
  -1 & -1  
  \end{array}
  \right).
\]
Hence this is different from the flux of the catenoid cousin.
\end{Exa}

\begin{Rmk}
  As discussed in Theorem 3.2 of 
  \cite{ruy1}, branched \cmcone{} immersions $x$ are divided 
  into the following   three classes: 

\medskip
\centerline{\rm $\bullet$\, irreducible,\quad
$\bullet$\,$\Hyp^1$-reducible,\quad $\bullet$ $\Hyp^3$-reducible.}

\medskip
If the surface $x$ is $\Hyp^3$-reducible,
  then the lift $F$ itself is single-valued, hence the integral 
  $({-1}/{2 \pi i}) \int F^{-1}\cdot dF$ on any homology class is 
  well-defined.  
  We call this flux $\Fl^\#$, 
  and it also satisfies a balancing formula on the ends. 
  (It is in fact, the $\Fl$ flux for $x^\#$.  Unlike the $\Fl$ flux, 
  $\Fl^{\#}$ does depend on the choice of Weierstrass data.)  
  In the case that $x$ is $\Hyp^1$-reducible, we cannot 
  define $\Fl^\#$, but there exists a special Weierstrass data $(g,\omega)$
  such that 
  $g\omega$ is single-valued, as seen in \cite{UY1}. So we can define a 
  flux $\fl^\#$ 
  on homology classes $\gamma$ by 
  $\fl^\#(\gamma)=(-1/{2 \pi i}) \int_{\gamma} g\,\omega.$ %
  When  $x$ is $\Hyp^3$-reducible, $\fl^\#$ is just one of the 
  diagonal components of $\Fl^\#$.
  The computation of these new fluxes at irregular ends is 
easier than that of the flux $\Fl$,
  because 
 $F^{-1}\cdot dF$ has at most a pole
whereas
$dF\cdot F^{-1}$ has an essential singularity.
\end{Rmk}

As an application of the flux formula,
we give the following result:  
\begin{Thm}\label{thm:1.8}
  Let $\overline{M}$ be a closed Riemann surface
  and consider a branched \cmcone{} immersion 
  $x:\overline{M}\setminus\{p_1,\ldots,p_n\}\to \Hyp^3$
  whose ends are all regular.
  Suppose that $p_1$ is a Type I end.
  Suppose also that the ends $\{p_1,\ldots,p_n\}$ are all poles or zeros of 
  the hyperbolic Gauss map $G$.
  Then $x$ has at least one more Type I end $z=p_j$\,\,$(2\leq j\leq n)$. 
\end{Thm}

Catenoid cousins and
4-ended Costa cousins constructed in \cite{ruy1} 
are examples satisfying 
the assumptions
of the theorem.
The assumption that the ends are all poles or zeros of the hyperbolic 
Gauss map is essential in this theorem. 
(See the following Example \ref{kumamotoexample}.)

\begin{Proof*}{Proof of Theorem~\ref{thm:1.8}}
Consider the 1-form
$
  \eta:=G\cdot \omega^\#,
$
where $\omega^\# := - Q/dG$.  
Since all ends of $x$ are regular, $\eta$ is a meromorphic 
1-form on $\overline{M}$.
Moreover, $\eta$ and $-\eta$ are 
the diagonal components of the $\sl(2,\C)$-valued
1-form $-dF\cdot F^{-1}$, so $\eta$ is holomorphic on 
$\overline{M}\setminus\{p_1,\ldots,p_n\}$.
Thus we have 
$
  \sum_{j=1}^n \eta_j=0,
$
where $\eta_j$ is the residue of $\eta$ at $p_j$.
Suppose the Hopf differential $Q$ has at most poles of order $1$
at $p_2,\ldots,p_n$.
Then $Q$ has the following Laurent expansions
at the ends
\begin{equation}
  Q(z)=\left(
         \frac{q_{-2}^{(j)}}{(z-p_j)^2}+\frac{q_{-1}^{(j)}}{z-p_j}+\cdots
         \right)\,dz^2,
\qquad (j=1,\ldots,n),
\end{equation}
where $q_{-2}^{(1)}\ne 0$ and $q_{-2}^{(j)}=0$ 
 for $j=2,\ldots,n$.
Since each end $z=p_j$\,\,($j=1,\ldots,n$) is a pole or a zero
of $G$,  we have
\[
  {G}/{dG}:=a_j(z-p_j)+b_j(z-p_j)^2+\cdots \qquad (j=1,\ldots,n)
\]
where $a_j$ is a non-zero complex number.
Thus we can conclude that
$\eta_1=-a_1q_{-2}^{(1)}(\ne 0)$ and
$\eta_j(z)=0$ for $j=2,\ldots,n$. 
This contradicts the relation $\sum_{j=1}^n \eta_j=0$.
\end{Proof*}

\begin{Exa}\label{kumamotoexample}
It can be easily seen that there is a complete \cmcone{} immersion
$x:\C \setminus \{0,1\}\to \Hyp^3$ such that
  \[ 
     G = \frac{z+1}{(z-1)^3} (z^2 -4z +1) \; , \quad g = z^2 \; , \quad
     Q = \frac{2}{z(z-1)^2}\, dz^2 \; , 
  \]
which has three regular ends and only $z=1$
is an end of Type I.  
 This shows that the condition $\{p_1 ,\ldots,p_n \}$ are 
        all poles or zeros of 
        $G$ is essential in Theorem~\ref{thm:1.8}.
Moreover, just like the flux $\Fl$, the KKMS-flux 
        {\em cannot\/} 
        always be nonzero on Type II ends.  This follows because 
        the KKMS-flux satisfies a balancing formula,
        and it can be easily verified that the KKMS-flux for regular
        ends of 
        Type I does not vanish.  
       (For minimal surfaces in $\R^3$, the KKS-flux 
             vanishes for ends of Type II.)  
        This suggests that for \cmcone{} surfaces, the KKMS-flux, 
        in addition to being difficult to compute,
        has no particular advantage over the flux $\Fl$. 
        It also lends support to the possibility of a relationship between 
        the two fluxes.  
\end{Exa}


\paragraph*{\bf  Further Remarks}
Let $\Met(\overline M)$ be the set of 
positive semi-definite conformal metrics of constant curvature $1$ 
with conical singularities on a closed Riemann surface $\overline M$.
Suppose that $du^2\in \Met(\overline M)$
has conical singularities at points $z_j \in \overline{M}$
($j=1,\ldots,n$) 
with order $\beta_j$ ($>-1$), that is, 
the metric admits a tangent cone of angle $2\pi(\beta_j+1)>0$
at each $z_j$. 
Following Troyanov \cite{Troyanov1}, we call 
$\beta=\sum_{j=1}^n \beta_j z_j$
the divisor represented by $du^2$.
As shown in \cite{UY3}, 
there is a canonical correspondence
between pseudo-metrics in $\Met({\overline M})$ and
(not totally umbilic) branched \cmcone{} surfaces 
in  $\Hyp^3$
with prescribed hyperbolic Gauss map
on $\overline M$.
Our first motivation to formulate Theorem~\ref{thm:1.8} was 
based on this problem:

\begin{Prob}
  Let $G$ be a meromorphic function on a closed 
  Riemann surface $\overline{M}$ with the ramification divisor
  \[
    \beta_0=m_1p_1+\cdots+m_n p_n \qquad (n,m_1,...,m_n\in \Z^{+}),
  \]
  where $p_1,\dots,p_n\in \overline{M}$ 
  are branch points of $G$. 
  Can one show the non-existence of conformal metrics
  $d\sigma^2_t\in \Met(\overline{M})$ 
  with the divisor 
  \[
    \beta_t=(t+m_1)p_1+\cdots+m_n p_n 
  \qquad \qquad \mbox{for $t\not\in\Z$?}
\]
\end{Prob}
We cannot remove the assumption $t\not\in {\bf Z}$.
In fact, the two pull-backs of the Fubini-Study metric $G^*d\sigma^2_0$
and $g^*d\sigma^2_0$ in Example \ref{kumamotoexample} have the divisors $1p_0+2p_1+1p_\infty$
and $1p_0+1p_\infty$, respectively.
When $\overline{M}=P^1(\C)$, a much stronger assertion is true: 
there are no metrics $d\sigma^2_t$ 
with the divisor $\beta_t$ even when $\beta_0$ is not the 
ramification divisor of some meromorphic function
(see, for example, Appendix of \cite{ruy2}).
However, on a closed surface $\overline{M}$ of genus$\ge 1$,
there exists a metric in $\Met(\overline{M})$ 
with only one conical singularity, by a result in \cite{Troyanov2}.  
This problem seems to be related with flux formulas for minimal 
surfaces in
$\R^3$ and for CMC-$1$ surfaces in ${\cal H}^3$ by the following two 
observations:

\begin{Obs}
  There is no smooth deformation of 
  metrics $du^2_t\in \Met(\overline{M})$ for $t\in\R$
  close to $0$ such that
  \begin{enumerate}
    \item $du^2_0=G^*d\sigma_0^2$, where $d\sigma^2_0$ is the
          Fubini-Study metric on $P^1(C)$.
    \item Each $du^2_t$ has the divisor $\beta_t$.
  \end{enumerate}
\end{Obs}

\begin{Obs}
  If the assertion of Theorem~\ref{thm:1.8} would have held 
  without assuming that
  $\{p_1,\dots,p_n\}$ are all poles or zeros of $G$, then
  there would be no such metrics as in the problem above.
\end{Obs}

\noindent
{\it Proof.}
  As seen in [UY3;p211], the difference of Schwarzian derivatives
  \[
     Q_t(z):=\frac{\widetilde S(du^2_t)-\widetilde S(du^2_0)}{2t}
  =\frac12 \left(\int_0^1\frac{\partial h_{ts}(z)}{\partial t}ds\right)dz^2
\]
  is a well-defined, holomorphic $2$-differential on 
  $\overline{M}\setminus\{p_1,\dots,p_n\}$, where $\tilde S(du_r^2)=h_r(z)dz^2$
for $r\in \R$.
   Then by \cite{Troyanov1}, 
the limit $Q_0:=\lim_{t\to 0}Q_t$
  has a pole of order $2$ 
  at $p_1$ and could possibly have poles of order at most one 
  at $p_2 ,\dots,p_n$.
  On the other hand, there exists a unique branched \cmcone{} surface 
  $x_t:\overline{M}\setminus\{p_1,\dots,p_n\}\to \Hyp^3(-t^2)$ with
  the hyperbolic Gauss map $G$ and with the associated metric
  $du^2_t$ for each $t$ (see \cite[Theorem 2.2]{UY3}).
  We identify the hyperbolic space ${\cal H}^3(-t^2)$ 
  of constant curvature $-t^2$ with the Poincar\'e ball
  \[
    \left(
      \left\{x\in\R^3\,|\,|x|<2/t\right\}\, , \,
      \frac{4}{4-t^2|x|^2}\,|dx|^2
    \right).
  \]
  Then $x_0=\lim_{t\to0}x_t:M^2\to \R^3$ is a conformal
  branched immersion with Hopf differential $Q_0$ and 
  Gauss map $G$. (See [UY2].)   $p_1$ is the only
  Type I end of the immersion $x_t$\,\, ($t\in \R$).
This makes a contradiction when $t=0$, since
the other ends of $x_0$
  have vanishing flux, since a minimal surface has vanishing flux 
  at Type II ends.
This proves Observation 1.
Next we suppose 
the assertion of Theorem \ref{thm:1.8} holds without assuming that
  $\{p_1,\dots,p_n\}$ are all poles or zeros of $G$.
Then we have the same contradiction for $x_t$ for each 
fixed $t\not\in \Z$, which proves the second observation. 

Unfortunately, as seen in Example~\ref{kumamotoexample}, 
Theorem~\ref{thm:1.8}
fails to hold without assuming that
$\{p_1,\dots,p_n\}$ are all poles or zeros of $G$.
However, the two observations suggest that there are 
connections between the set $\Met(\overline{M})$ and
the flux formula on branched conformal CMC-$1$ immersions.
Several interesting relations 
between these two objects are discussed in [RUY2] and [UY6].

\par
\setlength{\parindent}{0in}
\footnotesize
\vspace{0.2cm}
{\sc Wayne Rossman}:\,Graduate School of Mathematics,
Kyushu University,
Fukuoka 812-81, Japan.\,{\it E-mail}: 
wayne@math.kyushu-u.ac.jp\par
\vspace{3mm}
{\sc Masaaki Umehara}:\,Department of Mathematics,
Graduate School of Science,
Osaka University,
Toyonaka, Osaka 560,
Japan.\,
{\it E-mail}:
umehara@math.wani.osaka-u.ac.jp
\par
\vspace{3mm}
{\sc Kotaro Yamada}:\,
Department of Mathematics,
Faculty of Science,
Kumamoto University,
Kumamoto 860,
Japan.\,
{\it E-mail}: 
kotaro@gpo.kumamoto-u.ac.jp

\footnotesize
\begin{thebibliography}{KKMS}

\bibitem[Bie]{Bie}
  L.~Bieberbach,
  {\sc  Theorie der gew\"ohnlichen 
    Differentialgleichungen},
  Zweite Auflage.
  Springer-Verlag, 1965.
\bibitem[Bry]{Bry}
  R.~Bryant,
  {\it Surfaces of mean curvature one in hyperbolic space},
  Ast\'erisque {\bf 154--155} (1987), 341--347.
\bibitem[KKS]{kks} 
  N.~Korevaar, R.~Kusner, and B.~Solomon,
  {\it The structure of 
       complete embedded surfaces with constant mean curvature}, 
   J. Differential Geometry {\bf 30} (1989), 465--503.
\bibitem[KKMS]{kkms}
  N.~Korevaar, R.~Kusner, W.~H.~Meeks and B.~Solomon,
  {\it Surfaces in Hyperbolic space}, 
  Amer. J. Math. {\bf 114} (1992) 1--43.
\bibitem[RUY1]{ruy1}
  W.~Rossman, M.~Umehara and K.~Yamada,
  {\it Irreducible constant mean curvature 1 surfaces in
    hyperbolic space with positive genus}, 
  to appear in Tohoku Math. J.
\bibitem[RUY2]{ruy2}
  W.~Rossman, M.~Umehara and K.~Yamada,
  {\it Metrics with conical singularities and
  \cmcone{} surfaces of low total curvature
  in hyperbolic space}, 
  in preparation.
\bibitem[Sm]{Sm} 
  A.~J.~Small,
  {\it Surfaces of Constant Mean Curvature 
       $1$ in $H^3$ and Algebraic Curves on a Quadric}, 
  Proc.~Amer.~Math.~Soc. 
  {\bf 122} (1994), 1211--1220.
\bibitem[T1]{Troyanov1}
  M.~Troyanov,
  {\it Metric of constant curvature on a sphere with two conical
       singularities}, 
  in Diff. Geom., Lect. Notes in Math. vol.~1410, 
  Springer-Verlag,  (1989), 296--306. 
\bibitem[T2]{Troyanov2}
    M.~Troyanov,
    {\it Prescribing curvature on compact surfaces with conical
          singularities},
   Trans.~Amer.~Math.~Soc. {\bf 324} (1991), 793--821.
\bibitem[UY1]{UY1}
  M.~Umehara and K.~Yamada,
  {\it Complete surfaces of constant mean curvature-$1$
       in the hyperbolic $3$-space},
  {Ann. of Math. {\bf 137} (1993), 611--638.}
\bibitem[UY2]{UY2}
   M.~Umehara and K.~Yamada,
   {\it A parametrization of Weierstrass formulae and perturbation
        of some complete minimal surfaces of
        $\R^3$ into the hyperbolic $3$-space},
   J.~Reine Angew.~Math. {\bf 432} (1992), 93--116.
\bibitem[UY3]{UY3}
   M.~Umehara and K.~Yamada,
   {\it Surfaces of constant mean curvature-$c$
        in $H^3(-c^2)$ with prescribed hyperbolic Gauss map},
   Math. Ann. {\bf 304} (1996), 203--224.
\bibitem[UY4]{UY5}
   M.~Umehara and K.~Yamada,
   {\it Another construction of a \cmcone{} surface in $H^3$},
    Kyungpook Math. J. {\bf 35} (1996), 831--849.
\bibitem[UY5]{UY4}
    M.~Umehara and K.~Yamada,
   {\it A duality on \cmcone{} surface in the hyperbolic $3$-space
        and a hyperbolic analogue of the Osserman Inequality},
    Tsukuba J. Math. {\bf 21} (1997), 229-237.
\bibitem[UY6]{UY6}
    M.~Umehara and K.~Yamada,
   {\it Metric of constant curvature one with
        three conical singularities}, in preparation.
\bibitem[Y]{Yu}
    Z.~Yu,
   {\it Value distribution of hyperbolic Gauss maps},
    to appear in Proc.\ Amer.\ Math.\ Soc.
\end{thebibliography}
\end{document}